\documentclass{article}
\usepackage[T2A]{fontenc}
\usepackage[cp1251]{inputenc}
\usepackage[english, russian]{babel}
\usepackage[tbtags]{amsmath}
\usepackage{amsfonts,amssymb,amscd}
\usepackage{amsthm}
\usepackage[all]{xy}
\makeatletter
\@addtoreset{equation}{subsection}
\makeatother

\usepackage[OT2,T1]{fontenc}
\DeclareSymbolFont{cyrletters}{OT2}{wncyr}{m}{n}
\DeclareMathSymbol{\Sha}{\mathalpha}{cyrletters}{"58}

\newcommand{\comp}\circ

\renewcommand\labelenumii{\rm (\alph{enumii})}

\theoremstyle{plain}
\newtheorem{theorem}[subsection]{Теорема}
\newtheorem{lemma}[subsection]{Лемма}

\newtheorem{propos}[subsection]{Предложение}

\newtheorem{sled}[subsection]{Следствие}

\newtheorem*{claim*}{Claim}

\theoremstyle{definition}
\newtheorem{definition}[subsection]{Определение}

\newtheorem*{definition*}{Definition}

\newtheorem{example-remark}[subsection]{Remark-Example}
\newtheorem{subexample-remark}[equation]{Remark-Example}

\newtheorem*{notation*}{Notation}

\newtheorem{remark}[subsection]{Замечание}

\newcounter{NN}

\newcounter{NO}

\begin{document}

\title{Автоморфизмы трёхмерных многообразий, представимых в виде пересечения двух квадрик}
\author{А.\,А.~Авилов}

\maketitle

\begin{abstract}
Доказывается, что все трёхмерные $G$-многообразия дель Пеццо степе-\newline ни 4 с терминальными особенностями, за исключением однопараметрического семейства и четырёх выделенных случаев, эквивариантно перестраиваются в проективное пространство $\mathbb{P}^{3}$, квадрику $Q\subset \mathbb{P}^{4}$, $G$-расслоение на коники или поверхности дель Пеццо. Также мы покажем, что одно из четырёх выделенных многообразий является бирационально жёстким относительно подгруппы в группе автоморфизмов индекса 2.

Библиография: 14 названий.
\end{abstract}

\markright{Автоморфизмы трёхмерных пересечений двух квадрик}
\footnotetext[0]{Работа выполнена при частичной финансовой поддержке РФФИ (гранты № 15-01-02164 и 15-01-02158).}

\section{Введение}
Одной из мотивировок данной работы является получение бирациональной классификации трёхмерных $G$-многообразий, т.е. многообразий с бирегулярным действием конечной группы $G$. Мы будем работать над алгебраически замкнутым полем характеристики 0. Хорошо известно, что в этом случае существуют эквивариантное разрешение особенностей (см., например, ~\cite{5}) и эквивариантная программа минимальных моделей в размерности 3 (см., например,~\cite{8} и ~\cite[\S 3.6]{88}), которая позволяет привести любое трёхмерное $G$-многообразие $X$ с помощью определённых эквивариантных бирациональных перестроек к $G$-многообразию $\bar{X}$ со следующими свойствами: оно $G\mathbb{Q}$-фактори\-ально, имеет не более чем терминальные особенности и либо канонический класс является эффективным, либо на нём есть структура $G$-расслоения Мори.
\begin{definition} $G$-многообразие $\bar{X}$ с не более чем терминальными $G\mathbb{Q}$-фак\-ториальными особенностями называется \emph{$G$-расслоением Мори}, если существует такой морфизм $\pi:\bar{X}\to Y$, что $\pi_{*}\mathcal{O}_{\bar{X}}=\mathcal{O}_{Y}$, $\dim X>\dim Y$, $\rho^{G}(X/Y)=1$ и относительный антиканонический класс $-K_{X/Y}$ является $\pi$-обильным.
\end{definition}
\begin{definition} Пусть $\pi:\bar{X}\to Y$ -- $G$-расслоение Мори. Если $Y$ является точкой, то многообразие $\bar{X}$ называется \emph{$G\mathbb{Q}$-многообразием Фано}. Если при этом канонический класс является дивизором Картье, то $\bar{X}$ называется \emph{$G$-многообразием Фано}.
\end{definition}

\begin{definition} Пусть $\pi:\bar{X}\to Y$ -- трёхмерное $G$-расслоение Мори. Тогда оно называется \emph{расслоением на коники} (соотв., \emph{поверхности дель Пеццо степени $d$}), если общий слой расслоения изоморфен конике (соотв., поверхности дель Пеццо степени $d$).
\end{definition}

\begin{definition}\label{definition88} Трёхмерное многообразие $\bar{X}$ называется \emph{многообразием дель Пеццо} (см., например,~\cite{6}), если оно имеет не более чем терминальные горенштейновы особенности, а его антиканонический класс $-K_{\bar{X}}$ является обильным дивизором Картье и делится на 2 в группе Пикара. Если $G$ -- такая конечная подгруппа $\operatorname{Aut}(\bar{X})$, что $\bar{X}$ является  $G\mathbb{Q}$-многообразием Фано, то будем говорить, что $\bar{X}$ \emph{$G$-минимально}, а группу $G$ в этом случае будем называть \emph{минимальной}.
\end{definition}

\begin{definition} Трёхмерное многообразие $\bar{X}$ называется \emph{слабым многообразием дель Пеццо}, если оно имеет не более чем терминальные горенштейновы особенности, а его антиканонический класс $-K_{\bar{X}}$ является численно эффективным объёмным  дивизором Картье и делится на 2 в группе Пикара.
\end{definition}

Трёхмерные $G$-многообразия дель Пеццо были частично классифицированы Ю. Прохоровым в работе ~\cite{2}. Основным их инвариантом является \emph{степень} $d=(-\frac{1}{2}K_{\bar{X}})^{3}$, которая может принимать значения от 1 до 8. В этой работе мы рассмотрим случай $d=4$. Этот выбор обусловлен тем, что в случаях $d\geq 5$ имеется ровно четыре $G$-многообразия Фано, и их группы автоморфизмов хорошо изучпны, поэтому $d=4$ -- первый нетривиальный случай.

Второй мотивировкой данной работы является изучение конечных подгрупп в группе Кремоны $\operatorname{Cr}_{3}(\mathbb{K})$, где $\mathbb{K}$ -- алгебраически замкнутое поле характеристики 0. Группа $\operatorname{Cr}_{n}(\mathbb{K})$ -- это группа бирациональных автоморфизмов проективного пространства $\mathbb{P}^{n}$. Конечные подгруппы $\operatorname{Cr}_{2}(\mathbb{K})$ были полностью классифицированы И. Долгачёвым и В. Исковских в работе ~\cite{1}. Основная суть метода классификации состоит в следующем. Пусть $G$ -- конечная подгруппа в $\operatorname{Cr}_{2}(\mathbb{K})$. Тогда действие $G$ регуляризуется, т.е. существует гладкое проективное многообразие $Z$, на котором $G$ действует \emph{бирегулярными} автоморфизмами с эквивариантным бирациональным отображением $Z\dasharrow \mathbb{P}^{2}$. Применив далее эквивариантную программу минимальных моделей, мы получим $G$-расслоение Мори, которое является либо $G$-расслоением на коники над $\mathbb{P}^{1}$, либо $G$-минимальной поверхностью дель Пеццо. Классифицировав все возможные минимальные группы для расслоений на коники и для поверхностей дель Пеццо, Долгачёв и Исковских получили полную классификацию конечных подгрупп в $\operatorname{Cr}_{2}(\mathbb{K})$. Но довольно часто полученные подгруппы являются сопряжёнными в $\operatorname{Cr}_{2}(\mathbb{K})$, поэтому их естественно отождествить. Несложно видеть, что $G$-многообразия $Z_{1}$ и $Z_{2}$ дают сопряжённые подгруппы в том и только том случае, когда есть $G$-эквивариантное бирациональное отображение $Z_{1}\dasharrow Z_{2}$. Поэтому кроме классификации всех рациональных $G$-расслоений Мори необходимо исследовать также и бирациональные отображения между различными расслоениями.

Действуя таким же образом в трёхмерном случае, можно свести задачу классификации конечных подгрупп в $\operatorname{Cr}_{3}(\mathbb{K})$ к задаче описания всех рациональных $G\mathbb{Q}$-расслоений Мори и эквивариантных бирациональных отображений между ними. Эта программа была реализована в некоторых частных случаях: например, классифицированы простые неабелевы группы, вкладывающиеся в $\operatorname{Cr}_{3}(\mathbb{K})$ (~\cite{10}, см. также ~\cite{4}, ~\cite{14}), а также $p$-элементарные подгруппы $\operatorname{Cr}_{3}(\mathbb{K})$ (см. ~\cite{11}). Нас будет интересовать следующий вопрос: какими должны быть многообразие дель Пеццо степени 4 и конечная минимальная группа, действующая на нём, чтобы не было $G$-эквивариантного отображения на $\mathbb{P}^{3}$, квадрику $Q\subset \mathbb{P}^{4}$, расслоение на коники или поверхности дель Пеццо с \emph{регулярным} действием группы $G$? Все перечисленные классы многообразий являются более простыми с точки зрения классификации $G$-расслоений Мори, чем оставшиеся $G$-многообразия Фано, поэтому вопрос естественный.

В работе мы используем следующие обозначения для групп:
\begin{itemize}
\item
$C_{n}$ -- циклическая группа порядка $n$;
\item
$D_{2n}$ -- диэдральная группа порядка $2n$;
\item
$S_{n}$ -- симметрическая группа степени $n$;
\item
$G^{n}$ -- прямое произведение $n$ копий группы $G$.
\end{itemize}

Основными результатами данной работы являются следующие теоремы:
\begin{theorem}\label{th10} Пусть $X$ -- $G$-многообразие дель Пеццо степени 4. Предположим, что $X$ не является $G$-бирационально эквивалентным $\mathbb{P}^{3}$ и квадрике в $\mathbb{P}^{4}$ с регулярным действием группы $G$, а также $G$-расслоению Мори над базой положительной размерности. Тогда $X$ является одним из следующих многообразий:
\begin{enumerate}
\item
пересечение двух квадрик в $\mathbb{P}^{5}$ с $\operatorname{rk}\operatorname{Cl}(X)=5$. Такое многообразие единственно (см. ~\cite{2}), а его полная группа автоморфизмов изоморфна \newline$(\mathbb{C}^{*}\rtimes C_{2})^{3}\rtimes S_{3}$. Оно подробно описано в разделе 5;
\item
гладкое пересечение двух квадрик. В этом случае возможны следующие варианты:
\begin{enumerate}
\renewcommand\labelenumii{(\roman{enumii})}
\item
$\operatorname{Aut}(X)\simeq C_{2}^{5}\rtimes C_{5}$;
\item
$\operatorname{Aut}(X)\simeq C_{2}^{5}\rtimes D_{12}$;
\item
$\operatorname{Aut}(X)\simeq C_{2}^{5}\rtimes D_{6}$;
\item
группа $\operatorname{Aut}(X)$ вкладывается в точную последовательность $$0\to C_{2}^{5}\to \operatorname{Aut}(X)\to S_{4}\to 0.$$
\end{enumerate}
В случаях \textnormal{(2, i)}, \textnormal{(2, ii)} и \textnormal{(2, iv)} многообразие $X$ единственно с точностью до изоморфизма. В случае \textnormal{(2, iii)} такие многообразия $X$ образуют однопараметрическое семейство.
\end{enumerate}
\end{theorem}
\begin{remark} В случаях \textnormal{(2, i)}-\textnormal{(2, iv)} мы не утверждаем, что $G=\operatorname{Aut}(X)$.
\end{remark}
\begin{theorem} Пусть $X$ -- многообразие из пункта \textnormal{(2, i)} теоремы ~\ref{th10}, а $G\simeq C_{2}^{4}\rtimes C_{5}$. Тогда оно является $G$-бирационально жёстким, т.е. если есть другое $G$-расслоение Мори $X'\to Y'$ с $G$-эквивариантным бирациональным отображением $X\dasharrow X'$, то $X'\simeq X$. Как следствие, любое другое $G$-расслоение Мори даёт нам несопряжённое вложение $G\subset\operatorname{Cr}_{3}(\mathbb{K})$.
\end{theorem}
\begin{remark} В случаях \textnormal{(1)} и \textnormal{(2, ii)}-\textnormal{(2, iv)} теоремы~\ref{th10} не утверждается, что многообразие $X$ нельзя $G$-эквивариантно перестроить в $\mathbb{P}^{3}$, квадрику в $\mathbb{P}^{4}$ или $G$-расслоение Мори над базой положительной размерности. Вопрос о существовании таких перестроек для некоторых подгрупп $G\subset \operatorname{Aut}(X)$ остаётся открытым.
\end{remark}
\section{Пересечения двух квадрик и символы Сегре}
Пусть $X$ -- трёхмерное многообразие дель Пеццо степени 4 (см. определение~\ref{definition88}). Напомним, что мы предполагаем, что $X$ имеет только терминальные горенштейновы особенности. Хорошо известны следующие факты:
\begin{theorem} (см. ~\cite[Corollary 1.7]{6})\label{th11} Многообразие дель Пеццо $X$ степени 4 является пересечением двух квадрик в $\mathbb{P}^{5}$.
\end{theorem}
\begin{propos} (см., например, ~\cite[Example 10.3.1]{13}) Многообразие дель Пеццо $X$ степени 4 является рациональным.
\end{propos}
\begin{propos}\label{propos1} Многообразие дель Пеццо $X$ степени 4 является пересеченим двух гладких квадрик.
\end{propos}
\begin{proof} По теореме ~\ref{th11} многообразие $X$ является пересечением двух квадрик, обозначим их $Q_{1}$ и $Q_{2}$. По теореме Бертини общий элемент $Q$ пучка $\langle Q_{1}, Q_{2}\rangle$ неособ вне $X$. Так как $X=Q_{1}\cap Q_{2}$, то $Q$ может иметь особенности только в конечном множестве $\operatorname{Sing}(X)$. Поэтому общий элемент пучка квадрик может быть особым в том и только в том случае, когда $Q_{1}$ и $Q_{2}$ имеют общую особую точку. В таком случае $X$ является пересечением двух конусов с общей вершиной, поэтому размерность касательного пространства в вершине равна пяти. С другой стороны, размерность касательного пространства в терминальной горенштейновой особой точке на трёхмерном многообразии равна четырём (см. ~\cite[Theorem 1.1]{12}). Противоречие.
\end{proof}
Для описания групп автоморфизмов многообразий дель Пеццо степени 4 рассмотрим более общую ситуацию пересечения двух квадрик произвольной размерности.

Рассмотрим многообразие $X=Q_{1}\cap Q_{2}\subset \mathbb{P}^n$, где $Q_{1}$ и $Q_{2}$ -- различные квадрики, причём квадрика $Q_{2}$ неособа. Обозначим через $\mathcal{P}$ пучок квадрик $$\mathcal{P}=\{Q_{\lambda, \mu}=\lambda Q_{1}+\mu Q_{2}, (\lambda : \mu)\in \mathbb{P}^{1}\}$$ (будем обозначать квадрику и её уравнение, а также матрицу соответствующей квадратичной формы одним символом).
\begin{definition}\label{definition1}
\emph{Дискриминантом пучка квадрик $\mathcal{P}$} называется многочлен степени $n+1$ от двух переменных $$\Delta=\Delta (\lambda, \mu)=\operatorname{Det}(\lambda Q_{1}+\mu Q_{2})$$.
\end{definition}
Дискриминант пучка квадрик зависит от выбора порождающих $Q_{1}$ и $Q_{2}$, однако его корни (с учётом кратностей) определены однозначно, с точностью до автоморфизма $\mathcal{P}\simeq\mathbb{P}^{1}$.

Пусть $(\bar{\lambda}:\bar{\mu})$ -- корень уравнения $\Delta=0$. Существует такое целое число $d\geq 0$, что все миноры матрицы $Q_{\bar{\lambda}:\bar{\mu}}$ порядка $n+1-d$ зануляются, но не все миноры порядка $n-d$. Обозначим через $l_{i}, i=0, 1, ..., d$ минимальную кратность корня $(\bar{\lambda}:\bar{\mu})$ в минорах порядка $n+1-i$, тогда $l_{i}>l_{i+1}$. Пусть $e_{i}=l_{i}-l_{i+1}$, где $0\leq i\leq d-1$ и $e_{d}=l_{d}$.
\begin{definition}\label{definition2}
Числа $e_{i}$ называются \emph{характеристическими числами корня} $(\bar{\lambda}:\bar{\mu})$.
\end{definition}

Пусть $(\lambda_{i}:\mu_{i}), i=1, 2, ..., r$ -- все корни уравнения $\Delta=0$, $e^{i}_{j}, j=0, 1, ..., d_{i}$ -- их характеристические числа, причём если $i_{1}<i_{2}$, то $d_{i_{1}}\geq d_{i_{2}}$, а в случае равенства  наборы характеристических чисел упорядочены лексикографически.
\begin{definition} \label{definition3}\emph{Символом Сегре} пересечения двух квадрик $X$ (или пучка квадрик $\mathcal{P}$) называется набор чисел $$\sigma_{X}=\sigma_{P}=[(e^{1}_{0}...e^{1}_{d_{1}}), (e^{2}_{0}...e^{2}_{d_{2}}), ..., (e^{r}_{0}...e^{r}_{d_{r}})].$$
\end{definition}
\begin{remark} Будем опускать скобки в символе Сегре, если в них стоит ровно одно число.
\end{remark}
\begin{remark} Каждому корню дискриминанта (скобке в символе Сегре) соответствует особая квадрика, являющаяся конусом с $d$-мерной вершиной, где~$d$ -- количество характеристических чисел, соответствующих данному корню (будем называть это число длиной скобки).
\end{remark}

\begin{theorem} \label{th1}(~\cite[Chapter XIII, \S 10, Theorem I]{9}) Два пучка квадрик $\mathcal{P}_{1}$ и~$\mathcal{P}_{2}$ изоморфны тогда и только тогда, когда существует автоморфизм $\mathbb{P}^{1}$, переводящий корни $\Delta_{\mathcal{P}_{1}}$ в корни $\Delta_{\mathcal{P}_{2}}$, причём наборы характеристических чисел у соответствующих корней совпадают.
\end{theorem}

 Таким образом, пучок квадрик однозначно определяется конфигурацией корней уравнения $\Delta=0$ и символом Сегре. Поэтому можно определить нормальную форму пучка квадрик $\mathcal{P}$.

Для произвольного числа $e^{i}_{j}$ из символа Сегре многообразия $X$ рассмотрим две $e^{i}_{j}\times e^{i}_{j}$-матрицы
$$Q_{1, i, j}=\begin{pmatrix}
0&0& ... & 1 & -\frac{\mu_{i}}{\lambda_{i}}\\
0& ... &1&  -\frac{\mu_{i}}{\lambda_{i}}&0\\
...&...&...&...&...\\
1& -\frac{\mu_{i}}{\lambda_{i}}&0&...&0\\
 -\frac{\mu_{i}}{\lambda_{i}}&0&0&...&0
\end{pmatrix},\ Q_{2, i, j}=\begin{pmatrix}
0&0&...&0&1\\
0&...&0&1&0\\
...&...&...&...&...\\
0&1&0&...&0\\
1&0&0&...&0
\end{pmatrix}.$$
\begin{sled} \label{sled1}
Можно выбрать однородные координаты в $\mathbb{P}^5$ таким образом, что в них матрицы $Q_{1}$ и $Q_{2}$ имеют блочно-диагональный вид $$Q_{1}=\operatorname{Diag}(Q_{1, 1,1}, ..., Q_{1, r, d_{r}}), \  Q_{2}=\operatorname{Diag}(Q_{2, 1,1},..., Q_{2, r, d_{d}}).$$
\end{sled}

Теперь вернёмся к ситуации трёхмерного многообразия дель Пеццо $X$ степени 4. Мы знаем, что оно является пересечением двух квадрик, одна из которых неособа, поэтому ему можно сопоставить его символ Сегре.

\begin{propos} \label{propos2} Пусть $X$ является многообразием дель Пеццо степе-\newline ни 4. Тогда любая скобка в символе Сегре многообразия $X$ содержит не более двух характеристических чисел, а любая скобка из двух чисел имеет вид $(a, 1)$.
\end{propos}
\begin{proof}
Действительно, если в какой-то скобке содержится более двух чисел, то квадрика, соответствующая корню дискриминанта с этим набором характеристических чисел, является конусом над коникой с двумерной вершиной. Это вершина пересекается с другой квадрикой из пучка по кривой особых точек, что противоречит терминальности $X$. Любой скобке вида $(a, b)$, соответствует конус над неособой квадратичной поверхностью с одномерной вершиной. Простая проверка (см. следствие ~\ref{sled1}) показывает, что если $b>1$, то эта вершина целиком лежит на $X$, что даёт прямую особых точек.
\end{proof}
\begin{remark} Если же все скобки в символе Сегре многообразия $X$ имеют вид $(a)$ или $(a, 1)$, то несложно проверить, что особыми точками на $X$ будут вершины конусов, соответствующих скобкам вида $(a),\ a>1$, и точки пересечения одномерных вершин конусов, соответствующих скобкам длины 2, с другой квадрикой из пучка (для скобок вида $(1, 1)$ пересечение состоит из двух точек, для скобок вида $(a, 1),\ a>1$ -- из одной). В частности, пересечение двух квадрик неособо тогда и только тогда, когда его символ Сегре равен $[1, 1, 1, 1, 1, 1]$.
\end{remark}
\begin{remark} Кроме того, существует ровно одно многообразие с символом Сегре $[(1,1), (1,1), (1,1)]$. Несложно проверить (например, написав явные уравнения, см. следствие ~\ref{sled1}), что оно совпадает с многообразием из пункта~1 теоремы ~\ref{th10}.
\end{remark}
\begin{theorem}\label{th2} Пусть $X$ -- трёхмерное многообразие дель Пеццо степени 4. Тогда группа $\operatorname{Aut}(X)$ вкладывается в точную последовательность $$0\to \operatorname{Aut}(X)'\to \operatorname{Aut}(X)\to \operatorname{Aut}(X)''\to 0,$$ где группа $\operatorname{Aut}(X)'$ действует на $\mathbb{P}^{5}$, сохраняя каждую квадрику из пучка, а $\operatorname{Aut}(X)''$ -- группа автоморфизмов $\mathbb{P}^{1}$, переводящая каждый корень дискриминанта в корень дискриминанта с тем же набором характеристических чисел.
\end{theorem}
\begin{proof} Это простое следствие из теоремы ~\ref{th1}.
\end{proof}

Далее мы рассмотрим многообразие дель Пеццо $X$ степени 4 с действием $G$ -- такой конечной подгруппы $\operatorname{Aut}(X)$, что многообразие $X$ является $G$-минималь\-ным.

\begin{theorem}\label{th3}Пусть $X$ -- $G$-многообразие дель Пеццо. Предположим, что $X$ не является $G$-бирационально эквивалентным $\mathbb{P}^{3}$, квадрике в $\mathbb{P}^{4}$ и $G$-эквивари\-антному расслоению Мори с базой положительной размерности. Тогда символ Сегре $X$ равен $[1, 1, 1, 1, 1, 1]$ или $[(1,1), (1,1), (1,1)]$.
\end{theorem}
\begin{proof} Согласно утверждению ~\ref{propos2}, символ Сегре многообразия~$X$ имеет только скобки вида $(a)$ или $(a, 1)$.
\begin{lemma}\label{lemma1}
Если в символе Сегре многообразия $X$ есть ровно одна скобка вида $(n), n>1$ (соотв., $(n, 1)$), то $X$ является $G$-бирационально эквивалентным квадрике в $\mathbb{P}^{4}$ (соотв., $G$-расслоению на коники).
\end{lemma}
\begin{proof}
Рассмотрим случай скобки вида $(n)$. Обозначим соответствующий ей конус через $Q_{1}$, а его вершину через $p$. Тогда точка $p$ является особой $G$-инвариантной точкой многообразия $X$. Рассмотрим проекцию из этой точки. Общая прямая, лежащая на $Q_{1}$ и проходящая через $p$, пересекается с другой квадрикой $Q_{2}$ (которую можно считать гладкой) в двух точках, одна из которых $p$. В противном случае, любая образующая конуса либо целиком лежит на $Q_{2}$, либо имеет пересечение кратности 2 в точке $p$. Таким образом, любая образующая конуса $Q_{1}$ лежит на касательной плоскости в точке $p$ к $Q_{2}$, чего быть не может. Таким образом, проекция из точки $p$ является $G$-эквивариантным бирациональным отображением на гиперповерхность степени 2 в $\mathbb{P}^{4}$.

Пусть теперь скобка имеет вид $(n, 1)$. Обозначим соответствующий ей конус через $Q_{1}$, а его вершину через $l$ (она является прямой). Рассмотрим проекцию из $l$. Образом этой проекции будет неособая квадрика $Q_{1}'\subset\mathbb{P}^{3}$ -- основание конуса $Q_{1}$. Общая плоскость, проходящая через $l$, пересекает $Q_{2}$ по неособой конике. Действительно, если бы сечение общей плоскостью имело особенность, то по теореме Бертини это была бы точка пересечения $l$ с $Q_{2}$. Тогда общее сечение было бы парой прямых, проходящих через $l\cap Q_{2}$, а это означает, что $Q_{2}$ -- конус с вершиной в $l\cap Q_{2}$. Противоречие. Таким образом, проекция из $l$ даёт нам структуру $G$-расслоения на коники над $Q_{1}'$. Разрешив его особенности и применив эквивариантную относительную программу минимальных моделей, мы получаем искомое $G$- расслоение Мори на коники.
\end{proof}
Таким образом, можно считать, что любая скобка кроме $(1)$ либо не входит в символ Сегре, либо входит более одного раза. Перечислим все возможные символы Сегре, удовлетворяющие этому свойству, учитывая, что сумма всех чисел в символе Сегре равна шести: $$[1, 1, 1, 1, 1, 1],\ [2, 2, 1, 1],\ [2, 2, 2],\ [3, 3],\ [(1, 1),(1, 1), 1, 1],$$ $$[(1, 1), (1, 1),(1,1)],\ [(2, 1), (2, 1)].$$

Если символ Сегре многообразия $X$ равен  $[2, 2, 1, 1], [3, 3]$ или $[(2, 1), (2, 1)]$, то $X$ содержит ровно две особые точки, причём прямая $l$, проходящая через них, содержится в $X$ (это делается явной проверкой, соответствующие уравнения описаны в следствии ~\ref{sled1}) и является $G$-инвариантной, а само $X$ является пересечением двух конусов с вершинами в этих точках. Обозначим эти конуса через $Q_{1}$ и $Q_{2}$. Проекция из прямой $l$ даёт бирациональное $\bar{G}$-отображение на~$\mathbb{P}^{3}$. Действительно, рассмотрим общую плоскость, содержащую $l$. Её пересечение с $Q_{i}$ равно $l+l_{i}$. Для общей плоскости прямые $l_{1}$ и $l_{2}$ пересекаются в одной точке, не лежащей на $l$, что и требовалось доказать.

Если символ Сегре равен $[2, 2, 2]$, то $X$ содержит ровно три особые точки, причём плоскость, проходящая через них, содержится в $X$ (это также делается явной проверкой с помощью следствия ~\ref{sled1}). Таким образом, многообразие~$X$ содержит $G$-инвариантную плоскость и не может быть $G$-минимальным.

Наконец, если символ Сегре равен $[(1, 1), (1, 1), 1, 1]$, то $X$ содержит 4 особые точки. Они не лежат на одной плоскости, что видно из уравнений $X$, см. следствие ~\ref{sled1}. Рассмотрим проекцию из трёхмерного проективного пространства, порождённого этими точками. Мы получаем $G$-эквивариантное расслоение над $\mathbb{P}^{1}$ на рациональные поверхности, являющиеся пересечиями двух квадрик. Применив эквивариантное разрешение особенностей расслоения, а затем эквивариантную относительную программу минимальных моделей, мы получим $G$-расслоение Мори на коники или поверхности дель Пеццо.
\end{proof}
\section{Факты об особенностях линейных систем}
Пусть $X$ -- трёхмерное многообразие с не более чем терминальными особенностями, $G$ -- конечная группа, действующая на $X$, причём $X$ является многообразием $G\mathbb{Q}$-Фано. Более того, мы считаем, что канонический дивизор $K_{X}$ является дивизором Картье (поскольку нам потребуются исключительно приложения к пересечениям двух квадрик). Пусть $X'\to Y'$ -- другое $G$-расслоение Мори. Предположим, что существует бирациональное $G$-эквивариантное отображение $f:X\dasharrow X'$.

Существует метод, позволяющий разложить $f$ в композицию элементарных отображений, называемых линками Саркисова (подробности см. в ~\cite{7}). Для этого выберем очень обильный $G$-инвариантный дивизор $M'$ на $X'$ и положим $\mathcal{M}=f^{-1}_{*}(|M'|)$. Ввиду того, что $X$ является многообразием $G\mathbb{Q}$-Фано, существует такое рациональное число $\mu$, что $\mathcal{M}\subset |-\mu K_{X}|$. Если $X'$ не изоморфно $ X$, то неравенства Нётера-Фано-Исковских (см. ~\cite[Theorem 2.4]{7}) дают нам неканоничность пары $(X, \frac{1}{\mu}\mathcal{M})$. Поэтому для доказательства бирациональной жёсткости $X$ достаточно описать все возможные неканонические центры пар вида $(X, \frac{1}{\mu}\mathcal{M})$, и для каждого неканонического центра описать соответствующий линк Саркисова. Если все полученные линки дают многообразие, изоморфное $X$, то оно является бирационально жёстким. Для описания нульмерных неканонических центров пар нам понадобится следующая теорема:
\begin{theorem}\label{th4}(см., например, ~\cite[Lemma 1.10]{7}) В этих условиях пусть точка $p\in X$ -- центр неканонической особенности пары $(X, \frac{1}{\mu}\mathcal{M})$. Пусть $Z=M_{1}\cdot M_{2}$ -- цикл, являющийся пересечением двух общих элементов линейной системы $\mathcal{M}$. Тогда $\operatorname{mult}_{p}Z>4\mu^{2}$.
\end{theorem}
Для изучения неканонических центров, являющихся кривыми, нам понадобится следующее утверждение:
\begin{propos}\label{th6} (см., например, ~\cite[Exercise 6.18]{3}) В наших условиях пусть неприводимая кривая $C$ является центром неканонической особенности пары $(X, \frac{1}{\mu}\mathcal{M})$. Тогда кратность $\mathcal{M}$ вдоль кривой $C$ больше $\mu$.
\end{propos}

В случае неканонического центра, являющегося гладкой кривой, для описания связанного с ней линка Саркисова необходима следующая теорема:

\begin{theorem}\label{th24}(см., например, ~\cite[Proposition 1.2]{15}) Пусть $Y$ и $X$ трёхмерные нормальные многообразия, а $f: E\subset Y\to C\subset X$ -- дивизориальное стягивание неприводимого дивизора $E$ на кривую $C$. Предположим, что $\dim f(Y^{\operatorname{sing}})=0$, многообразие $X$ имеет изолированные особенности, а $-E$ является $f$-обиль\-ным. Тогда $Y$ изоморфно раздутию $X$ вдоль $C$.
\end{theorem}

В частности, условия теоремы ~\ref{th24} выполнены в случае, когда многообразия $X$ и $Y$ имеют не более чем терминальные особенности, а $f$ является стягиванием Мори.
\section{Гладкие пересечения двух квадрик}

Пусть $X$ -- гладкое пересечение двух квадрик. Это равносильно тому, что символ Сегре $X$ равен $[1,1,1,1,1,1]$. В этом случае можно считать, что $$Q_{1}=\sum\limits_{i=1}^{6}\lambda_{i}x_{i}^{2},\ Q_{2}=\sum\limits_{i=1}^{6}x_{i}^{2},$$ где все $\lambda_{i}$ различны, а $x_{i},\ 1\leq i\leq 6$ -- некоторая система координат. Группа $\operatorname{Aut}(X)'$ (в обозначениях теоремы ~\ref{th2}) изоморфна $(C_{2})^5$ и действует обращением знаков у координат $x_{1}, ..., x_{5}$. Группа $\operatorname{Aut}(X)''$ является группой автоморфизмов $\mathbb{P}^{1}$, сохраняющих множество из шести точек $S=\{(\lambda_{i}:1)\mid i=1, ..., 6\}$. В общем случае эта группа тривиальна. Перечислим все случаи (с точностью до автоморфизма $\mathbb{P}^{1}$), в которых она нетривиальна: \begin{enumerate}
\renewcommand\labelenumii{(\roman{enumii})}
\item
$S=\{T_{0}T_{1}(T_{0}^{4}-T_{1}^{4})=0\}.$ В этом случае $\operatorname{Aut}(X)''\simeq S_{4}$.
\item
$S=\{T_{0}^6+T_{1}^{6}=0\}.$ В этом случае $\operatorname{Aut}(X)''\simeq D_{12}$.
\item
$S=\{T_{0}^6+aT_{0}^{3}T_{1}^{3}+T_{1}^{6}=0\}, a\neq -2, 0, 2.$ В этом случае $\operatorname{Aut}(X)''\simeq D_{6}$.
\item
$S=\{T_{0}T_{1}(T_{0}^4+aT_{0}^{2}T_{1}^{2}+T_{1}^{4})=0\}, a\neq -2, 0, 2.$ В этом случае $\operatorname{Aut}(X)''\simeq D_{4}$.
\item
$S=\{T_{0}(T_{0}^5+T_{1}^{5})=0\}.$ В этом случае $\operatorname{Aut}(X)''\simeq C_{5}$.
\item
$S=\{(T_{0}^2+T_{1}^{2})(T_{0}^2+aT_{1}^{2})(T_{0}^2+bT_{1}^{2})=0\}, a, b\neq -1, 0, 1, a\neq b.$ В этом случае $\operatorname{Aut}(X)''\simeq C_{2}$.
\end{enumerate}
Этот список можно получить из классификации конечных подгрупп $\operatorname{PGL}_{2}(\mathbb{C})$ и полуинвариантных бинарных форм относительно действия этих групп (см. ~\cite[\S 5.5]{1}) аналогично ~\cite[\S 6.4]{1}.

\begin{lemma}\label{le1} В случаях (4) и (6), а также в случае тривиальной группы $\operatorname{Aut}(X)''$, многообразие $X$ является $G$-бирационально эквивалентным $G$-рас\-слоению на поверхности дель Пеццо для любой подгруппы $G\subset\operatorname{Aut}(X)$.
\end{lemma}
\begin{proof} В этих случаях в $S$ есть одна, две или четыре орбиты, состоящие в совокупности из 4 точек, поэтому в $\mathbb{P}^{5}$ есть инвариантное трёхмерное подпространство, порождаемое вершинами соответствующих конусов (они имеют координаты $x_{i}=\delta_{i}^{j}$ для различных $j$). Проекция из него даёт искомое $G$-расслоение.
\end{proof}

Таким образом, для доказательства теоремы ~\ref{th10} в случае гладкого многообразия $X$ осталось проверить структуру группы $\operatorname{Aut}(X)$ в случаях (2), (3) и~(5). Зная уравнения многообразия $X$, в этих случаях легко явно указать образующие $\operatorname{Aut}(X)$  и проверить, что $\operatorname{Aut}(X)\simeq\operatorname{Aut}(X)'\rtimes \operatorname{Aut}(X)''$. В случае (1) не существует расщепляющего гомоморфизма $\operatorname{Aut}(X)''\to \operatorname{Aut}(X)$, поскольку можно явно проверить, что у элемента $\operatorname{Aut}(X)''$ порядка 4 не существует прообраза в $\operatorname{Aut}(X)$ порядка 4.

Докажем ещё следующую полезную лемму:
\begin{lemma}\label{lemma15} Пусть $X$ -- гладкое пересечение двух квадрик в $\mathbb{P}^{5}$. Пусть $G$ -- подгруппа $\operatorname{Aut}(X)$, имеющая неподвижную точку на $X$. Тогда $X$ является $G$-эквивалентным $G$-расслоению на квадрики.
\end{lemma}
\begin{proof} Пусть $\pi:\widetilde{X}\to X$ -- раздутие $X$ в неподвижной точке. Тогда~$\widetilde{X}$ является гладким $G$-многообразием дель Пеццо степени 3. Действительно, пусть $C$ -- произвольная кривая на $\widetilde{X}$. Если $C$ лежит на $E$, то $C\cdot K_{\widetilde{X}}<0$. В противном случае $$C\cdot K_{\widetilde{X}}=\pi(C)\cdot K_{X}+2C\cdot E=-2\deg \pi(C)+2\operatorname{mult}_{p}C\leq0.$$ Более того, равенство достигается только в случае, когда $\pi(C)$ -- прямая. Таким образом, численная эффективность  $-K_{\widetilde{X}}$ доказана, а его объёмность очевидна. Антиканоническая линейная система оботражает $\widetilde{X}$  на кубическую гиперповерхность в $\mathbb{P}^{4}$, содержащую $G$-инвариантную плоскость. Проекция из неё даёт искомую структуру $G$-расслоения на квадрики.
\end{proof}

В следующем разделе мы покажем, что в случае (5) пересечение двух квадрик является бирационально жёстким относительно группы $C_{2}^{4}\rtimes C_{5}$, а значит и относительно всей группы $\operatorname{Aut}(X)$. Остаётся открытым вопрос, есть ли в случаях (1)-(3) подгруппы в $\operatorname{Aut}(X)$, относительно которых $X$ является бирационально жёстким? Большую часть подгрупп $\operatorname{Aut}(X)$ можно отбросить, используя леммы \ref{le1} и \ref{lemma15}, но, к сожалению, не все.

\subsection{Случай $\operatorname{Aut}(X)''\simeq C_{5}$}
Разберём теперь подробнее случай (5). В этом случае $X$ можно задать системой уравнений $$\xi x_{1}^{2}+\xi^{2}x_{2}^{2}+\xi^{3}x_{3}^{2}+\xi^{4}x_{4}^{2}+x_{5}^{2}=\sum\limits_{i=1}^{6}x_{i}^{2}=0,$$ где $\xi$ -- корень пятой степени из 1. Группа $\operatorname{Aut}(X)$ изоморфна полупрямому произведению $(C_{2})^{5}\rtimes C_{5}\simeq C_{2}\times(C_{2}^{4}\rtimes C_{5})$, причём $(C_{2})^{5}$ действует сменой знаков у $x_{1},..., x_{5}$, а $C_{5}$ переставляет эти координаты по циклу.

\begin{propos}\label{propos3} Группа $\operatorname{Aut}(X)$ имеет следующие подгруппы:
$$\{e\},\ C_{2},\ C_{2}^{2},\ C_{2}^{3},\ C_{2}^{4},\ C_{2}^{5},\ C_{5},\ C_{10},\ C_{2}^{4}\rtimes C_{5},\ C_{2}^{5}\rtimes C_{5}.$$
\end{propos}
\begin{proof}Для произвольной подгруппы $G\subset \operatorname{Aut}(X)$ положим $G'=G\cap\operatorname{Aut}(X)'$ и $G''=\pi(G)$, где $\pi:\operatorname{Aut}(X)\to \operatorname{Aut}(X)''$ -- стандартная проекция (в обозначениях теоремы \ref{th2}). Если $G''$ тривиальна, то $G=G'\simeq(C_{2})^n$.

Предположим теперь, что $G''=C_{5}$. Рассмотрим произвольный элемент нашей подгруппы порядка 5. Несложно показать, что с помощью сопряжения элементом из $(C_{2})^{5}$ его можно перевести в элемент, переставляющий координаты по циклу без смены знаков. Отождествим группу $(C_{2})^{5}$ с векторным пространством $\mathbb{F}_{2}^{5}$ над полем $\mathbb{F}_{2}$, на этом пространстве задано тавтологическое представление группы $C_{5}$. Подгруппа $G'$ соответствует подпредставлению, которых ровно 4: нульмерное, тривиальное одномерное (порожденное вектором $(1, 1, 1, 1, 1)$), ортогональное ему четырёхмерное и пятимерное. Следовательно, с точностью до сопряжения $\operatorname{Aut}(X)$ содержит четыре таких подгруппы: $C_{5}$, $C_{5}\times C_{2}$, $(C_{2})^{4}\rtimes C_{5}$, $G$. Более того, подгруппа, изоморфная $C_{2}^{4}\rtimes C_{5}$, в группе $G$ ровно одна.
\end{proof}

\begin{propos} Пусть $X$ -- гладкое пересечение двух квадрик в $\mathbb{P}^{5}$ с группой автоморфизмов $\operatorname{Aut}(X)\simeq C_{2}^{5}\rtimes C_{5}$. Пусть $G$ -- подгруппа $\operatorname{Aut}(X)$, изоморфная $C_{5}$ или $C_{10}$. Тогда $X$ является $G$-бирационально эквивалентным $G$-расслоению на квадрики.
\end{propos}
\begin{proof} Это немедленно следует из леммы \ref{lemma15}.
\end{proof}

Обозначим подгруппу $C_{2}^{4}\rtimes C_{5}\subset \operatorname{Aut}(X)$ через $G$, а нормальную подгруппу $C_{2}^{4}\subset \operatorname{Aut}(X)$ через $G'$.

\begin{lemma}\label{lemma57} Пусть $Y\subset X$ -- сечение $X$ гиперплоскостью $x_{6}=0$, а $Z\subset X$ -- сечение $X$ гиперплоскостью $x_{i}=0,\ i\neq 6$. Тогда $G$-орбита точки на $Y$ может иметь длину 16, 20, 40 или 80, а $G'$-орбита точки на $Z$ состоит из 4, 8 или 16 точек.
\end{lemma}
\begin{proof} Пусть $y$ -- некоторая точка $Y$, а $G_{y}$ -- её стабилизатор. Если в $G_{y}$ есть элемент порядка 5, то первые пять координат $y$ ненулевые, поэтому никакой элемент $G'$ в стабилизаторе лежать не может, и орбита $y$ имеет длину 16. Если $G_{y}$ не содержит элементов порядка 5, то $G_{y}\subset G'$. Поскольку среди координат $y$ не менее трёх ненулевых, то мощность $G_{y}$ равна 1, 2 или 4. Вторая часть утверждения доказывается аналогично.
\end{proof}

\begin{theorem} Пусть $X$ -- гладкое пересечение двух квадрик в $\mathbb{P}^{5}$ с группой автоморфизмов $\operatorname{Aut}(X)\simeq C_{2}^{5}\rtimes C_{5}$. Пусть $G\simeq C_{2}^{4}\rtimes C_{5}$ -- подгруппа $\operatorname{Aut}(X)$. Тогда $X$ является $G$-бирационально жёстким.
\end{theorem}
\begin{proof} Пусть $\mathcal{H}\subset |-\mu K_{X}|$ -- некоторая $G$-инвариантная линейная система без неподвижных компонент. Докажем сначала, что точка не пожет быть неканоническим центром пары $(X, \frac{1}{\mu}\mathcal{H})$. Предположим противное -- пусть точка
$$p=(p_{1}:p_{2}:p_{3}:p_{4}:p_{5}:p_{6})\in X$$
является неканоническим центром. Из уравнений многообразия $X$ легко выводится, что среди чисел $p_{i}$ как минимум три не равны нулю. Без ограничения общности можно считать, что $p_{1}\neq0$ и $p_{2}\neq 0$. Тогда точки
$$(\pm p_{1}:\pm p_{2}:\pm p_{3}:p_{4}:p_{5}:p_{6}),$$ где количество минусов чётно, лежат в $G$-орбите $p$, поэтому тоже являются неканоническими центрами. Все эти 4 точки лежат на двумерной плоскости $S$, которая не имеет других точек пересечения с $X$. Рассмотрим общую гиперплоскость $H$, содержащую $S$, и два общих элемента линейной системы $H_{1}, H_{2}\in\mathcal{H}$. Пусть $Z=H_{1}\cdot H_{2}$. Из того, что $S\cap X$ нульмерно следует, что $H$ не содержит компонент $Z$. Тогда
$$16\mu^{2}=H\cdot Z\geq4\operatorname{mult}_{p}Z>16\mu^{2},$$ где последнее неравенство следует из теоремы ~\ref{th4}. Полученное противоречие показывает, что точка не может быть неканоническим центром.

Пусть $L=L_{1}$ -- неприводимая кривая, являющаяся неканоническим центром, $d=\deg L_{1}$ -- её степень, а $\{L_{i}\mid 1\leq i\leq r\}$ -- её $G$-орбита. Пусть $H$ -- достаточно общее гиперплоское сечение $X$, а $H_{1}$ и $H_{2}$ -- общие элементы $\mathcal{H}$. Тогда
$$16\mu^{2}=H\cdot H_{1}\cdot H_{2}\geq (\operatorname{mult}_{L}\mathcal{H})^{2}H\cdot L>\mu^{2}rd.$$ Из этого следует, что $rd\leq 15$. Поскольку $r$ является индексом некоторой подгруппы $G$, то (см. утверждение ~\ref{propos3}) $r$ может быть равным 1, 5 или 10.

Предположим, что $r=5$. В таком случае кривая $L_{1}$ сохраняется группой~$G'$, а её степень не превосходит 3. Выберем такое $1\leq i\leq 5$, что $L_{1}$ не лежит в гиперплоскости $\{x_{i}=0\}$. Тогда пересечение $L_{1}$ с этой гиперплоскостью является множеством из не более чем трёх точек, инвариантным относительно~$G'$. Согласно лемме \ref{lemma57}, такого не может быть. Противоречие.

Пусть $r=10$. В этом случае $d=1$. Без ограничения общности можно считать, что $L_{1}$ и $L_{2}$ образуют $G'$-орбиту. Найдётся такое $i$, что $L_{1}$ и $L_{2}$ не лежат на гиперплоскости $\bar{H}=\{x_{i}=0\}$, иначе $L_{1}$ или $L_{2}$ лежала бы на плоскости $\{x_{i}=x_{j}=x_{k}=0\}$ для некоторых различных $i, j, k$, но пересечение этой плоскости с $X$ состоит из 4 точек. Тогда $(L_{1}\cup L_{2})\cap \bar{H}$ является парой точек, инвариантной относительно $G'$, чего не может быть.

Остался последний случай $r=1$. В этом случае $L$ не лежит в гиперплоскости $\{x_{1}=0\}$. Пересечение $L$ с гиперплоскостью $\{x_{1}=0\}$ является $G'$-инвариантным множеством, поэтому из леммы \ref{lemma57} следует, что степень кривой может быть равна 4, 8 или 12. Кроме того, $L$ лежит в гиперплоскости $\bar{H}=\{x_{6}=0\}$, поскольку иначе их пересечение было бы $G$-инвариантным конечным подмножеством $\bar{H}$, но минимальная мощность орбиты в этом случае равна $16>d$ (см. лемму \ref{lemma57}). Легко понять, что в подпространстве размерности меньше 4 кривая $L$ не лежит, поскольку соответствующее пятимерное представление $G$ неприводимо. Поэтому кривая $L$ не может иметь степень 4.

Кривая $L$ является $G$-инвариантной кривой на поверхности дель Пеццо $\bar{H}$. Пусть $M$ и $M_{i}, 1\leq i\leq 5$ -- стандартный базис группы Пикара поверхности $\bar{H}$, где $M_{i}$ -- непересекающиеся $(-1)$-кривые, а $M$ -- класс собственного прообраза прямой в $\mathbb{P}^{2}$ (напомним, что любая поверхность дель Пеццо степени 4 является раздутием $\mathbb{P}^{2}$ в пяти точках). Пусть $N$ -- любая $(-1)$-кривая на $\bar{H}$, тогда $L\cdot N$ не зависит от $N$, поскольку $L$ является $G$-инвариантной, а $G'$ действует на множестве $(-1)$-кривых транзитивно. Действительно, пусть элемент $g\in G'$ сохраняет некоторую $(-1)$-кривую, которая в нашем случае является прямой в $\mathbb{P}^{4}$. Элемент $g$ меняет знаки у двух или четырёх координат, во втором случае его неподвижные точки лежат на сечении $\bar{H}$ координатной гиперплоскостью, в первом -- на сечении $\bar{H}$ координатным подпространством коразмерности 2. Действие $g$ на прямой имеет две неподвижные точки, поэтому прямая обязана целиком лежать на сечении $\bar{H}$ координатной гиперплоскостью, чего не может быть, поскольку это сечение является неприводимой кривой степени 4.

 Поскольку индекс пересечения со всеми $(-1)$-кривыми один и тот же, то $$L=aM-b\sum M_{i}$$ для некоторого натурального $a$ и неотрицательного целого $b$. Пересечём теперь~$L$ с $(-1)$-кривой с классом $N=M-M_{1}-M_{2}$ и с гиперплоским сечением (его класс равен $-K_{\bar{H}}=3M-\sum M_{i}$). Получаем $$b=a-2b,\ 3a-5b=d.$$ Откуда выводим, что $$a=3b,\ b=\frac{d}{4}.$$ Таким образом, $L=-\frac{d}{4}K_{\bar{H}}$. Докажем, что $L$ является полным пересечением~$\bar{H}$ с гиперповерхностью в $\mathbb{P}^{4}$ степени $\frac{d}{4}$.

Любой дивизор из линейной системы $|-K_{\bar{H}}|$ является гиперплоским сечением, поскольку наше вложение $\bar{H}\subset \mathbb{P}^{4}$ является каноническим. Применив формулу Римана-Роха и теорему Кодаиры о занулении к дивизорам $-K_{\bar{H}}$, $-2K_{\bar{H}}$ и $-3K_{\bar{H}}$, получаем
$$5=h^{0}(-K_{\bar{H}})=\frac{-K_{\bar{H}}\cdot(-K_{\bar{H}}-K_{\bar{H}})}{2}+\chi(\mathcal{O})=4+\chi(\mathcal{O}),$$ $$h^{0}(-2K_{\bar{H}})=\frac{-2K_{\bar{H}}\cdot(-2K_{\bar{H}}-K_{\bar{H}})}{2}+\chi(\mathcal{O})=12+\chi(\mathcal{O}),$$ $$h^{0}(-3K_{\bar{H}})=\frac{-3K_{\bar{H}}\cdot(-3K_{\bar{H}}-K_{\bar{H}})}{2}+\chi(\mathcal{O})=24+\chi(\mathcal{O}).$$ Отсюда получаем, что $$h^{0}(-2K_{\bar{H}})=13,\ h^{0}(-3K_{\bar{H}})=25,$$ что совпадает с размерностями пространств сечений $\bar{H}$  квадратичными и кубическими гиперповерхностями соответственно.

Мы получили, что $L=\bar{H}\cap \{S=0\}$, где $S$ -- квадратичный или кубический многочлен. Напомним, что $\bar{H}$ задано уравнениями $$Q_{1}=\xi x_{1}^{2}+\xi^{2}x_{2}^{2}+\xi^{3}x_{3}^{2}+\xi^{4}x_{4}^{2}+x_{5}^{2}=0,\ Q_{2}=\sum\limits_{i=1}^{5}x_{i}^{2}=0.$$  Если многочлен $S$ кубический, то, прибавив к $S$ многочлен вида $Q_{1}P_{1}+Q_{2}P_{2}$ для некоторых линейных многочленов $P_{i}$, можно добиться того, что коэффициенты при мономах $x_{i}^{2}x_{i\pm1},\ 1\leq i\leq 6$ (мы отождествляем $x_{0}$ с $x_{6}$ и $x_{7}$ с~$x_{1}$) обнулятся, причём многочлены $P_{i}$ единственны (это несложно проверить, написав явные уравнения на неизвестные коэффициенты многочленов~$P_{i}$). Полученное уравнение является $G$-полуинвариантным, чего не бывает для кубических многочленов.

Пусть многочлен $S$ квадратичный. В этом случае можно считать, что $S$ полуинвариантен относительно действия группы $G$ (возможно, после прибавления к нему $p_{1}Q_{1}+p_{2}Q_{2}$ для некоторых $p_{i}$). В таком случае легко проверить, что $S$ имеет вид $$\zeta x_{1}^{2}+\zeta^{2}x_{2}^{2}+\zeta^{3}x_{3}^{2}+\zeta^{4}x_{4}^{2}+x_{5}^{2},$$ где $\zeta$ -- корень пятой степени из единицы, отличный от $\xi$. В частности, кривая~$L$ неособа. По теореме \ref{th24}, соответствующий линк Саркисова происходит из раздутия кривой $L$. Обозначим это раздутие через $f:\widetilde{X}\to X$, через $\widetilde{H}$ обозначим собственный прообраз $\bar{H}$, а через $E$ -- исключительный дивизор. На конусе Мори $\operatorname{NE}(\widetilde{X})$ есть два экстремальных луча, стягивание первого даёт морфизм~$f$, обозначим второй луч через $R$. Рассмотрим произвольную кривую $C$ на $\widetilde{H}$. Тогда $$C\cdot K_{\widetilde{X}}=C\cdot (f^{*}K_{X}+E)=f(C)\cdot K_{X}+C\cdot E=-2\deg C+f(C)\cdot L=0,$$ где последний индекс пересечения рассматривается на $\bar{H}$. Аналогично показывается, что для произвольной кривой $C\subset \widetilde{X}$ индекс пересечения $C\cdot K_{\widetilde{X}}\leq 0$. Таким образом, $K_{\widetilde{X}}\cdot R=0$, но стягивание $R$ малым не является, поэтому раздутие $L$ не даёт линка Саркисова (см. ~\cite{7}). Таким образом, линейная система~$\mathcal{H}$ не имеет неканонических центров, и теорема доказана.
\end{proof}

\begin{remark} Группа $G$ может действовать на других расслоениях Мори, например на $Y\times\mathbb{P}^{1}$, где $Y$ -- поверхность дель Пеццо степени 4 с группой автоморфизмов $C_{2}^{4}\rtimes D_{10}$ (см. ~\cite{1}). Таким образом, группа $\operatorname{Cr}_{3}(\mathbb{K})$ содержит как минимум две несопряжённые подгруппы, изоморфные $G$.
\end{remark}

\section{Пересечения двух квадрик с максимальным рангом $\operatorname{Cl}(X)$}
Пусть $X$ -- пересечение двух квадрик с символом $[(1,1),(1,1),(1,1)]$. Многообразие такого типа единственно. Можно считать, что $X$ задано уравнениями $$Q_{1}=x_{1}x_{2}+\xi x_{3}x_{4}+\xi^{2}x_{5}x_{6}=0,\ Q_{2}=x_{1}x_{2}+x_{3}x_{4}+x_{5}x_{6}=0,$$ где $\xi$ -- кубический корень из 1.

Множество особых точек $X$ состоит из шести обыкновенных двойных точек $\{x_{i}=\delta_{i}^{j}\}$ для $j=1,...,6$, $r=5$, группа $\operatorname{Cl}(X)$ порождается плоскостями, лежащими на $X$ (всего их 8, а их уравнения имеют вид $x_{i}=x_{j}=x_{k}=0$, где $i\in \{1, 2\},\ j\in \{3, 4\},\ k\in \{5, 6\}$), ранг группы $\operatorname{Cl}X$ равен пяти, более того, это единственное пересечение двух квадрик с $\operatorname{rk}\operatorname{Cl}X\geq5$ (см. ~\cite{2}).

Рассмотрим точную последовательность $$0 \longrightarrow \operatorname{Aut}^0(X) \longrightarrow\operatorname{Aut}(X) \longrightarrow \widehat{\operatorname{Aut}}(X) \longrightarrow 0,$$ где $\operatorname{Aut}^0(X)$  -- ядро действия на $\operatorname{Cl}(X)$, а группа $\widehat{\operatorname{Aut}}(X)$ действует эффективно на $\operatorname{Cl}(X)$.
Группа $\operatorname{Aut}^0(X)$ действует тривиально на множестве плоскостей, а следовательно и на множестве особых точек. Из этого нетрудно вывести, что она является трёхмерным тором. Согласно ~\cite[Corollary 7.4, Corollary 7.5 (i)]{2}, группа $\widehat{\operatorname{Aut}}(X)$ содержится в $W(\Delta'')\simeq C_{2}\times S_{4}$ (где $W(\Delta'')$ -- группа Вейля некоторой системы корней, канонически связанной с $X$, подробности см. в~\cite{2}). С другой стороны, группа $\widehat{\operatorname{Aut}}(X)$ содержит элементы, меняющие местами координаты $x_{2i-1}$ и $x_{2i},\ 1\leq i\leq 3$, и элементы $$g_{1}:(x_{1}:x_{2}:x_{3}:x_{4}:x_{5}:x_{6})\mapsto(x_{3}:x_{4}:x_{5}: x_{6}:x_{1}:x_{2}),$$ $$g_{2}:(x_{1}:x_{2}:x_{3}:x_{4}:x_{5}:x_{6})\mapsto(x_{1}:x_{2}:\xi x_{5}:\xi x_{6}:\xi^{2}x_{3}:\xi^{2} x_{4}),$$  которые в совокупности порождают подгруппу в $\operatorname{Aut}(X)$, изоморфную \newline $C_{2}^{3}\rtimes S_{3}\simeq C_{2}\times S_{4}$. Таким образом, $\operatorname{Aut}(X)\simeq (\mathbb{C}^{*})^{3}\rtimes (C_{2}^{3}\rtimes S_{3})$. Для нас будет удобно рассматривать  $\widehat{\operatorname{Aut}}(X)$ как подгруппу в $S_{6}$ (группе перестановок $x_{i}$) со следующими порождающими: $$h_{1}=(1324),\ h_{2}=(12),\ h_{3}=(56),\ h_{4}=(135)(246).$$

\begin{theorem}\label{th17} Пусть многообразие $X$ является $G$-минимальным для некоторой конечной подгруппы $\operatorname{Aut}(X)$. Тогда $G$ вкладывается в следующую точную последовательность: $$0\to C_{k}\times C_{l}\times C_{m}\to G\to G'\to 0,$$ где $G'\subset\widehat{\operatorname{Aut}}(X)$ -- одна из следующих групп: $$C_{4},\ C_{2}^{2},\ C_{4}\times C_{2},\ D_{8},\ C_{2}^{3},\ D_{8}\times C_{2},\ S_{4},\ C_{2}^{3}\rtimes C_{3},\ C_{2}^{3}\rtimes S_{3},$$ а $k,\ l,\ m$ -- натуральные числа (возможно, равные единице). Во всех случаях, кроме $D_{8}$, группа $G'$ единственна с точностью до сопряжения в $\widehat{\operatorname{Aut}}(X)$, есть ровно три класса сопряжённости в случае $G'\simeq D_{8}$.
\end{theorem}
\begin{proof}
Рассмотрим действие $G$ на множестве плоскостей. Это множество либо состоит из одной орбиты, либо разбивается на две орбиты, состоящие из плоскостей, лежащих в гиперплоскостях $x_{2i-1}=0$ и $x_{2i}=0$ соответственно (поскольку сумма плоскостей в каждой орбите должна быть пропорциональна каноническому классу), без ограничения общности можно ситать, что $i=3$.

Рассмотрим отображение $\pi: G\to \widehat{\operatorname{Aut}}(X)$, являющееся ограничением отображения $\operatorname{Aut}(X)\to \widehat{\operatorname{Aut}}(X)$, и обозначим через $G'$ его образ. В ядре отображения $\pi$ находится конечная подгруппа $(\mathbb{C}^{*})^{3}$, поэтому ядро изоморфно \newline $C_{k}\times C_{l}\times C_{m}$. Группа $G'$ также минимальна. Классифицируем все возможные минимальные подгруппы $\widehat{\operatorname{Aut}}(X)$ с точностью до сопряжения.

Сначала предположим, что мы в ситуации, когда множество плоскостей разбивается на две орбиты. Тогда $G'$ лежит в $D_{8}$ -- подгруппе $\langle g_{1}, g_{2}\rangle \subset\widehat{\operatorname{Aut}}(X)$. Минимальная подгруппа должна иметь порядок 4 или 8. Явной проверкой можно убедиться, что минимальными являются следующие подгруппы: $$A_{1}=\langle g_{1}\rangle,\ A_{2}=\langle g_{1}^{2},\ g_{2}\rangle,\ A_{3}=\langle g_{1}, g_{2}\rangle.$$

Теперь рассмотрим случай, когда все плоскости лежат в одной орбите, но порядок $G'$ не делится на 3. Тогда без ограничения общности можно считать, что $G'\subset \langle g_{1}, g_{2}, g_{3}\rangle \simeq D_{8}\times C_{2}$ (так как все силовские 2-подгруппы сопряжены этой, а нас интересуют подгруппы с точностью до сопряжения). В этом случае $G'$ является одной из следующих групп: $$A_{4}=\langle g_{1}, g_{3}\rangle,\ A_{5}=\langle g_{1}^{2}, g_{2}, g_{3}\rangle,\ A_{6}=\langle g_{1}, g_{2}, g_{3}\rangle,$$ $$A_{7}=\langle g_{1}g_{3}, g_{2}\rangle,\ A_{8}=\langle g_{1}, g_{2}g_{3}\rangle.$$

Теперь рассмотрим случай, когда все плоскости лежат в одной орбите, и порядок $G'$ делится на 3. Тогда либо $G'$ совпадает с $\widehat{\operatorname{Aut}}(X)$, либо является подгруппой индекса 2. Среди трёх подгрупп индекса 2 ровно две являются минимальными: $$A_{9}=\langle g_{1}, g_{2}g_{3}, g_{4}\rangle,\ A_{10}=\langle g_{1}^{2}, g_{2}, g_{3}, g_{4}\rangle.$$
\end{proof}
\begin{propos} Если в обозначениях теоремы \ref{th17} группа $G'$ не содержит элемента третьего порядка, то многообразие $X$ является $G$-эквивалентным $G$-расслоению Мори на коники или поверхности дель Пеццо.
\end{propos}
\begin{proof} В этом случае $G'$ сопряжена подгруппе $A_{i}$ для некоторого $i<9$. Из явного описания этих групп с помощью образующих видно, что множество из четырёх точек $$(1:0:0:0:0:0), (0:1:0:0:0:0), (0:0:1:0:0:0), (0:0:0:1:0:0)$$ является $G$-инвариантным. Следовательно, проекция из подпространства $\{x_{5}=x_{6}=0\}$ задаёт $G$-эквивариантное расслоение на рациональные поверхности. Применив эквивариантные разрешение особенностей и относительную программу минимальных моделей, получим искомое расслоение.
\end{proof}

Таким образом, если многообразие $X$ является $G$-бирационально жёстким, то либо $G'$ совпадает с $\widehat{\operatorname{Aut}}(X)$, либо сопряжена $A_{9}$ или $A_{10}$, то есть $G$ принадлежит одному из трёх семейств подгрупп.
\begin{remark} В этих случаях можно более точно описать группу $G''=\operatorname{ker}(G\to G')$: для некоторого $m$ имеется вложение $C_{m}^{3}\subset G''\subset C_{2m}^{3}$, причём $G''$ инвариантна относительно естественного действия группы $C_{3}$. Для каждого $m$ таких подгрупп ровно четыре. Это несложно получить из того, что $G''$ должна быть инвариантна относительно действия $G'$ сопряжениями.
\end{remark}
Остаётся открытым вопрос, для каких именно подгрупп $\operatorname{Aut}(X)$ многообразие $X$ является бирационально жёстким. Например, является ли оно бирационально жёстким относительно всей группы $\operatorname{Aut}(X)$?

\end{document}